\documentclass[a4paper]{article}
\usepackage[utf8]{inputenc}
\usepackage[]{graphicx}
\usepackage{hyperref}

\usepackage{amsmath}
\usepackage{amssymb}
\usepackage[margin=2.5cm,top=2.5cm,bottom=2.5cm]{geometry}

\usepackage{array}
\newcolumntype{L}[1]{>{\raggedright\let\newline\\\arraybackslash\hspace{0pt}}p{#1}}

\title{Rigorous computer-assisted bounds on renormalisation fixed point functions, eigenfunctions, and universal constants}
\author{Andrew Burbanks, Andrew Osbaldestin, Judi Thurlby }
\date{10 March 2021}

\begin{document}

\maketitle

\begin{abstract}
We gain tight rigorous bounds on the renormalisation fixed point function for period doubling in families of unimodal maps with degree 2 critical point. By writing the relevant eigenproblems in a modified nonlinear form, we use these bounds, together with a contraction mapping argument, to gain tight bounds on the essential eigenvalues and eigenfunctions of the linearised renormalisation operator at the fixed point and also those of the operator encoding the universal scaling of added uncorrelated noise.

We gain bounds on the corresponding power series coefficients and universal constants accurate to over 400 significant figures, confirming and (in the case of noise) extending the accuracy of previous numerical estimates, by using multi-precision interval arithmetic with rigorous directed rounding to implement operations on a space of analytic functions.
\end{abstract}

\section{Introduction}
Existence of the fixed point of the doubling operator for maps with degree 2 critical point was first proved by computer-assisted means by Lanford~\cite{Lan82}, by analytic means by Campanino and Epstein~\cite{Cam81} and, most recently, in full generalisation by Lyubich~\cite{Lyu99}. It is some decades since the first computer-assisted proofs. We revisit the problem with improved processing power, high precision computations, and  parallel processing, and extend its application to bound eigenfunctions of the linearised operator and also the operator controlling the scaling of uncorrelated noise.

In this paper we follow the method established in our paper~\cite{Bur20} where we studied the universality class corresponding to degree~4 critical points.  In section~\ref{sec:renormfp}, we work with a modified renormalisation operator, on a suitable space of analytic functions, corresponding to the action of the usual doubling operator on even maps. We identify a ball of functions centred on an approximate fixed point, given by a polynomial of degree $1280$, with $\ell^1$-radius $\rho\simeq 10^{-409}$.  We prove that a variant of Newton's method for the fixed-point problem (corresponding to our modified operator) is a contraction map on this ball, yielding rigorous  bounds on the fixed-point function itself and thereby on the associated universal constant, $\alpha$, controlling scaling in the state variable for families of maps in the corresponding universality class.

In section~\ref{sec:DTG}, we consider the eigenproblem for the derivative of our renormalisation operator at the fixed point.  We approach the eigenproblem in a novel way, by rewriting it in a modified nonlinear form, and apply a contraction mapping argument, using the ball of functions proven to contain the renormalisation fixed point, to bound the eigenfunctions corresponding to essential eigenvalues.  In particular, we gain tight bounds on the eigenfunction-eigenvalue pair for eigenvalue $\delta$, the universal constant controlling scaling in the parameter for the relevant families of unimodal maps.

In section~\ref{sec:LW}, we further adapt the technique to the eigenproblem corresponding to the universal scaling of added uncorrelated noise, gaining rigorous bounds on the relevant eigenfunction and hence on the eigenvalue $\gamma$.

Numerical approximations to the Feigenbaum constants $\alpha$ and $\delta$ (for the universality class of maps with degree $2$ critical points) with over $1,000$ digits have been computed by Broadhurst in 1999~\cite{Bro99}, and $10,000$ digits of each by Molteni in 2016~\cite{Mol16} using Chebychev polynomials. An approximation to $\gamma$ with $15$ digits was provided by Kuznetsov and Osbaldestin in 2002~\cite{Kuz02}. The rigorous bounds that we compute provide over $400$ confirmed digits for each constant (and for each one of $640$ nonzero coefficients of the relevant power series, together with bounds on all higher-order terms).

\section{The renormalisation fixed point}\label{sec:renormfp}

We consider the operator $R$ defined by: 
\begin{equation}
Rg(x):= a^{-1}g(g(ax)),
\end{equation}
where $a:= g(1)$ is chosen to preserve the normalisation $g(0)=1$. 

We first seek a nontrivial fixed point, $g^{*}$, of $R$, with a critical point of degree 2 at the origin.
It suffices to restrict to even functions. To this end, we define $X=Q(x)=x^2$,
and write
\begin{equation}
    g(x)=G(Q(x))=G(X),
\end{equation}
with $G$ in the Banach algebra $\mathcal{A}(\Omega)$ of functions analytic on an open disc $\Omega=D(c,r) :=\{z\subset \mathbb{C}:|z-c|<r\}$ and continuous on its closure, $\overline{\Omega}$, with finite $\ell^1$-norm. Specifically, we write $f\in\mathcal{A}(\Omega)$ as
\begin{equation}
f(z)=\sum_{k=0}^{\infty}a_k\left(\frac{z-c}{r}\right)^k,
\label{eqn:powerseries}
\end{equation}
i.e., we take the monomials $e_k:z\mapsto\left(\frac{z-c}{r}\right)^k$ as Schauder basis, with the corresponding $\ell^1$-norm,
\begin{equation}
\|f\|:=\sum_{k=0}^{\infty}|a_k|.\label{eqn:norm}
\end{equation}
We first define a modified operator $T$, corresponding to the action of $R$ on even functions, by 
\begin{equation}
TG(X):= a^{-1}G(Q(G(Q(a)X))),
\label{eq:TGX}
\end{equation}
in which $a:= G(1)$.
The above formulation is in contrast to~\cite{Lan82} in which the ansatz $g(x)=1+x^2h(x^2)$ is taken, with $h$ varying in a suitable space of functions equipped with a modified norm.

\begin{figure}[ht!]
\begin{center}
\includegraphics[width=7.5cm]{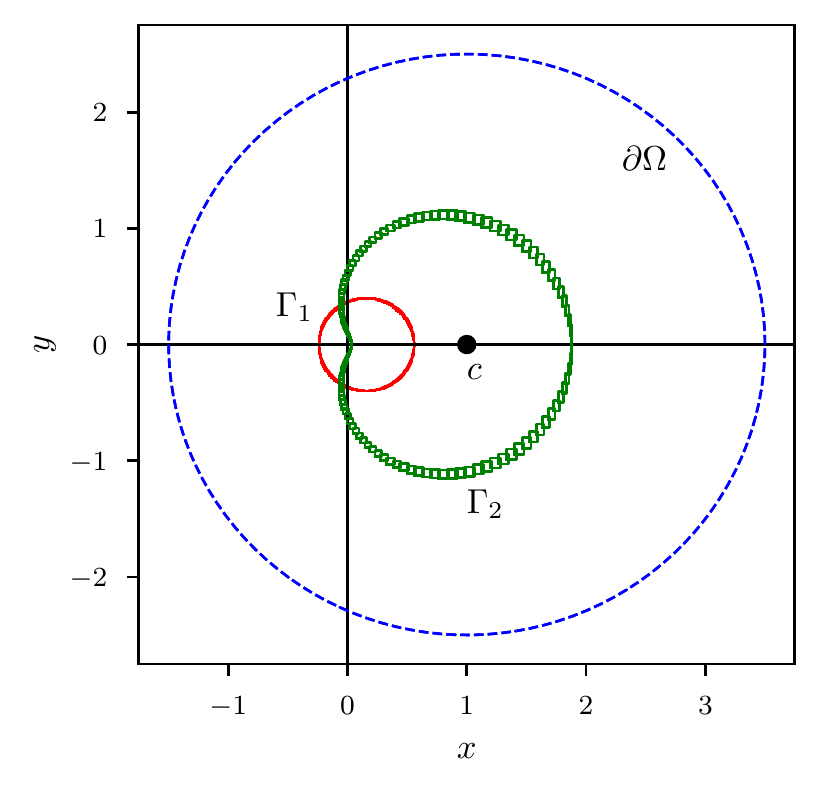}
\caption{Domain extension for the modified operator $T$ on $\mathcal{A}(\Omega)$ for the ball $B=B(G^{0},\rho)$, illustrated using a covering of the boundary $\partial\Omega\subset\mathbb{C}$ (dashed lines) of $\Omega=D(c,r)$ by $256$ rectangles, showing the corresponding coverings $\Gamma_1$ of $a^2\partial\Omega$ (red in colour copy) and $\Gamma_2$ of $Q(G(a^2\partial\Omega))$ (green in colour copy) valid for all $G\in B$.\label{fig:de}}
\end{center}
\end{figure}

We start by finding an accurate polynomial approximation, $G^0$, to the fixed point of $T$.
It is important, for what follows, to prove: (i) that the operator $T$ is well-defined on a certain $\ell^1$-ball, $B\subset\mathcal{A}(\Omega)$, of radius $\rho$ around $G^0$, (ii) that it is differentiable there, and (iii) that the derivative is compact.
We take domain $\Omega=D(1,2.5)$ to define the space $\mathcal{A}(\Omega)$, and establish that the `domain extension' or `analyticity-improving' property~\cite{Mac93} holds for our operator on the chosen ball in this space, i.e., that for all $G\in B$:
\begin{align}
\overline{a^2\Omega} &\subset \Omega,\label{eqn:de1}\\
\overline{Q(G(a^2\Omega))}&\subset\Omega,
\end{align}
in which the overline indicates topological closure.  (Note that, since $a=G(1)$, the universal quantifier on $G$ is not vacuous for equation~(\ref{eqn:de1}).)
This proves~\cite{Mac93} that the operator is well-defined on the ball $B$, that it is differentiable there, and that the derivative $DT(G)$ is compact for all $G\in B$. Figure~\ref{fig:de} demonstrates domain extension via a rigorous covering of the relevant sets.

We then use a contraction mapping argument on the ball $B$ to bound a fixed point of $T$. Note, however, that $T$ is not itself contractive at the fixed point we seek.
Following~\cite{Lan82} and~\cite{Eck84}, we  therefore consider a Newton-like operator corresponding to our modified operator $T$, 
\begin{equation}
    \Phi : G\mapsto G-\Lambda F(G),\label{eqn:newtonlike}
\end{equation}
where $F(G):=T(G)-G$ and $\Lambda$ is a fixed linear operator approximating  $[DF(G^0)]^{-1}$. We establish invertibility of $\Lambda$ and note that $\Phi$ therefore has the same fixed points as $T$.

To establish (uniform) contractivity of $\Phi$ on $B$, we bound a suitable norm of the derivative of $\Phi$:
\begin{equation}
    \|D\Phi(G)\|\leq\kappa<1\quad
    \text{for all }  G\in B,\label{eqn:contractive}
\end{equation}
and appeal to the mean value theorem.
Further, to establish that $\Phi$ is indeed a contraction map on $B$, we bound the movement of the approximate fixed point under  $\Phi$:
\begin{equation}
    \|\Phi(G^0)-G^0\|\leq\varepsilon,
\end{equation}
and establish the cruical inequality
\begin{equation}
    \varepsilon<\rho(1-\kappa),
\end{equation}
so that $\Phi(B)\subset B$.
The contraction mapping theorem then yields the existence of a (locally unique) fixed point $G^{*}\in B$, of $\Phi$, and therefore also of $T$.

In~(\ref{eqn:newtonlike}), above, we choose $\Lambda$ to be a fixed linear operator. The Fr\'echet derivative of $\Phi$ is thus given by
\begin{equation}
    D\Phi(G): \delta G \mapsto \delta G - \Lambda (DT(G)\delta G - \delta G).
    \label{eqn:DPhiG}
\end{equation}
Note that the Fr\'echet derivative, $DT(G)$, of our modified operator $T$ at $G$ is given by:
\begin{align}
    DT(G):\delta G\mapsto &-a^{-2}\delta a G(Q(G(a^2X)))\label{eq:term1}\\
    &+a^{-1}\delta G(Q(G(a^2X)))\label{eq:term2}\\
    &+a^{-1}G'(Q(G(a^2X)))\cdot 2G(a^2X)\cdot \delta G(a^2X)\label{eq:term3}\\
    &+a^{-1}G'(Q(G(a^2X)))\cdot 2G(a^2X)\cdot G'(a^2X)\cdot 2XG(1)\delta a\label{eq:term4},
\end{align}
where $\delta a=\delta G(1)$.
A simpler linear operator in which  terms corresponding to (\ref{eq:term1}) and (\ref{eq:term4}), which involve variations in $a$, are absent is often used in the literature when estimating the spectrum of $DR(g)$ nonrigorously. (In numerical calculations for the spectrum the simpler operator suffices with minor alterations to the spectral characteristics, see for instance \cite{Var11}.)

Note also that, in~(\ref{eqn:contractive}), we bound $D\Phi(G)$ by considering the `maximum column-sum norm': we note that
\begin{equation}
\|D\Phi(G)\|
:=
\sup_{\|f\|=1}\|D\Phi(G)f\|
\leq \sup_{k\ge 0}\|D\Phi(G)e_k\|,
\label{eqn:mcsnorm}
\end{equation}
where the norm on the left hand side of the inequality is the standard operator norm, and recall that $e_k$ denotes the $k$-th basis element.

For the rigorous calculations we use interval arithmetic, with rigorous directed rounding modes to bound operations in the corresponding space of analytic functions $\mathcal{A}(\Omega)$. The first detailed exposition of such a framework applied to renormalisation operators was provided in \cite{Eck84}, where it was applied with standard precision arithmetic to operations on functions of 2 real variables in the study of area-preserving maps.

For our computations, we introduce multi-precision arithmetic with directed rounding, specialise our implementation to spaces $\mathcal{A}(\Omega)$ and the corresponding complex Banach algebra, and use parallel computation to establish the bounds on contractivity. (For complex values, a straightforward analogue of interval arithmetic, namely rectangle arithmetic with intervals bounding real and imaginary parts, is used.)  Specifically, individual functions are written as the sum of a polynomial part (to some chosen truncation degree, $N$), and a high-order part, $f=f_P+f_H$. We maintain bounds on the power series coefficients such that those of $f_P$ lie in computer-representable intervals, and those of $f_H$ are bounded in norm, $\|f_H\|\le v_H$.

Following~\cite{Eck84}, in order to accommodate balls of functions and to absorb errors in the case where it would be undesirable to do so in the polynomial and high-order bounds, we write $f=f_P+f_H+f_E$, with $f_E$ a general `error' function bounded in norm, $\|f_E\|\le v_E$. We ensure that all computed operations deliver intervals bounding polynomial coefficients and upper bounds on the norms of the respective high-order and error parts that guarantee inclusion of the exact result.

The challenge of bounding the supremum in equation~(\ref{eqn:mcsnorm}) is reduced to a finite computation in two parts: Firstly, we bound $\|D\Phi(G)e_k\|$ (for all $G\in B$) for $k=0,1,\ldots,N$, by bounding the expressions in~(\ref{eqn:DPhiG})--(\ref{eq:term4}) evaluated at the polynomial basis elements. This first computation is well-suited to a parallel implementation, in which care is taken to ensure the safety of directed rounding modes across processes.  Secondly, we bound the action of $D\Phi(G)$ (for all $G\in B$) on a single ball $B_H$ of high-order functions $f_H$, such that $\|f_H\|\le 1$, that therefore contains all of the high-order basis elements $e_k$ for $k>N$. The latter requires careful consideration of the action of $D\Phi(G)$ on high-order perturbations $\delta G$ in order to minimise dependencies on $\delta G$ when implementing equation~(\ref{eqn:DPhiG}) in order to gain a suitable bound $\kappa<1$~\cite{Bur20}. (We additionally make use of closures in order to avoid recomputation of bounds on those subexpressions in the Fr\'echet derivative $D\Phi(G)\delta G$ that do not depend on $\delta G$.)
\begin{figure}[ht]
\begin{center}
\begin{tabular}{cc}
\includegraphics[width=7.5cm]{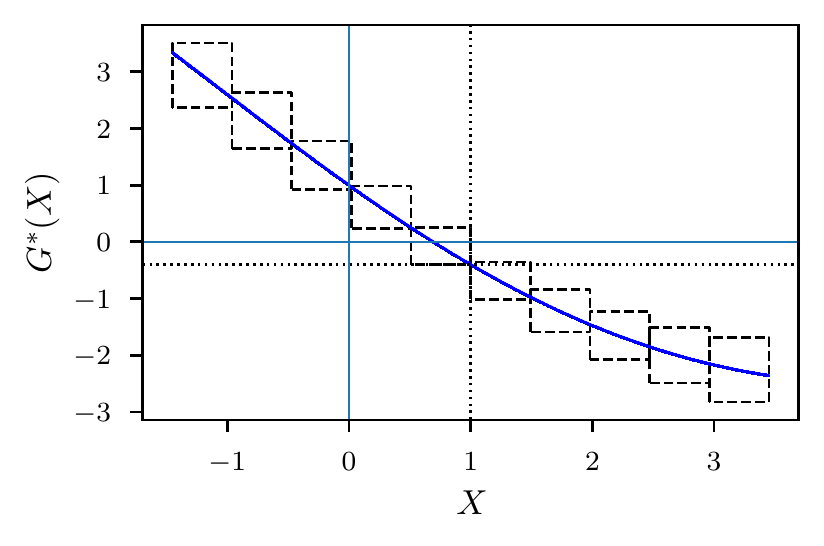}
& \includegraphics[width=7.5cm]{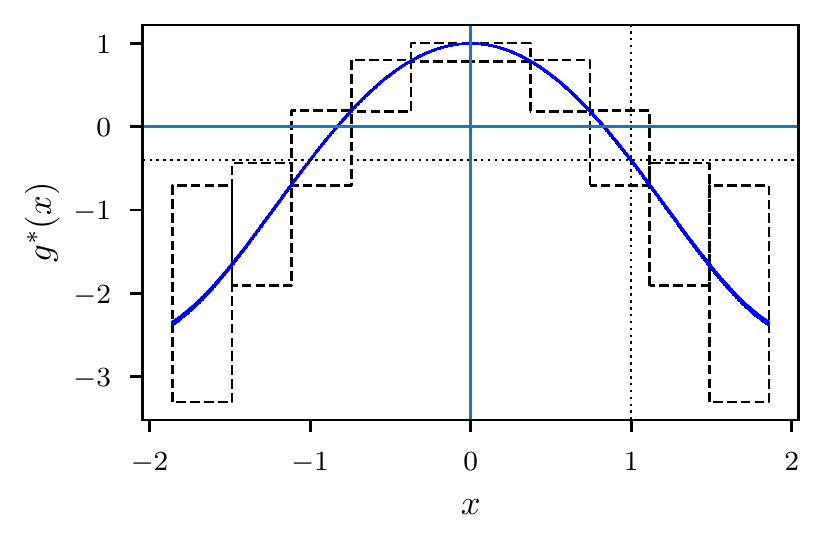}\\
(a) & (b)\\
\includegraphics[width=7.5cm]{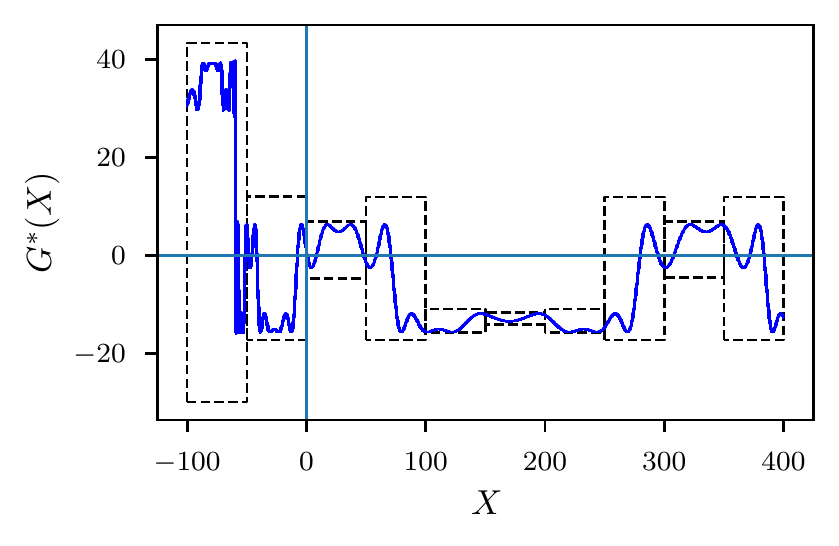}
& \includegraphics[width=7.5cm]{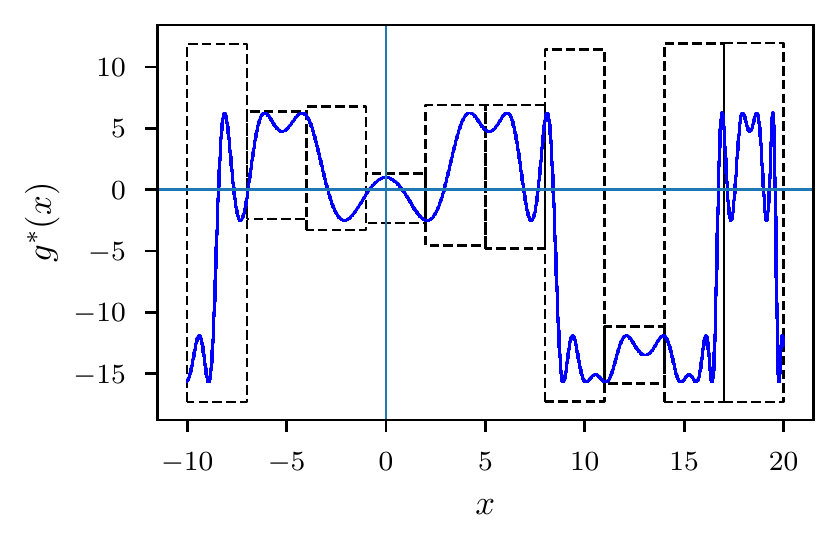}\\
(c) & (d)
\end{tabular}
\caption{Rigorous covering by rectangles of the graphs of the fixed point functions, (a,c) $G^{*}(X)$ of the operator $T$ and (b,d) $g^{*}(x)$ of the operator $R$. In (a) and (b) we use coverings of the interval $\Omega\cap\mathbb{R}$, respectively its preimage (for $X\ge 0$) under $Q$, by $10$ (dashed lines) and $1,000$ (solid lines) subintervals. The dotted lines indicate the universal constant $G^{*}(1)=g^{*}(1)=a\simeq -0.3995$. In (c) and (d) the functions $G^{*}$ and $g^{*}$ are extended to larger domains by making use of the fixed-point equations, $T(G^{*})=G^{*}$ and $R(g^{*})=g^{*}$, recursively.\label{fig:bigandlittleg}}
\end{center}
\end{figure}

It would be sufficient, in order to gain a proof of existence of the fixed point by this technique, to use truncation degree $N=20$ for the polynomial part of $G$ (thus degree $2N=40$ for $g$) with standard double-precision accuracy and careful use of directed rounding modes.  However, we tighten the resulting bounds significantly by increasing the truncation degree and using multi-precision interval arithmetic.

Using truncation degree $N=640$ for the polynomial part of $G$ (thus degree $2N=1280$ for $g$) and precision corresponding to (approximately) $P=\lfloor 2N/3\rfloor=426$ digits in the significand, we are able to choose a function ball of radius $\rho=10^{-409}$ (actually a close dyadic rational approximation to $10^{-409}$, with $1415$ bits in the significand) that results in bounds $\varepsilon<7\times 10^{-410}$, and $\kappa<1.3\times 10^{-99}$ (the numbers listed here are safely upwards-rounded decimal conversions of the corresponding dyadic rational  computer-representable bounds).  We confirm rigorously that $\varepsilon<\rho(1-\kappa)$, establishing that the Newton-like operator $\Phi$ corresponding to our modified operator $T$ is a contraction map on $B(G^0,\rho)$. The fixed-point functions $G^{*}$ and $g^{*}$ are illustrated in figure~\ref{fig:bigandlittleg}.

As an immediate consequence of computing such tight bounds on the fixed point,  we gain rigorous bounds on the universal constants $a=G^{*}(1)$ and $\alpha=a^{-1}$ to over $400$ significant digits each (see appendix \ref{d2summary}) with the first 20 given as:
\begin{align}
    a=G^{*}(1)&=-0.39953\,52805\,23134\,48985...\\
    \alpha= a^{-1}&=-2.50290\,78750\,95892\,8222...
\end{align}

\section{Eigenfunctions of the linearised operator at the fixed point}
\label{sec:DTG}

\begin{figure}[ht!]
\begin{center}
\begin{tabular}{cc}
\includegraphics[width=7.5cm]{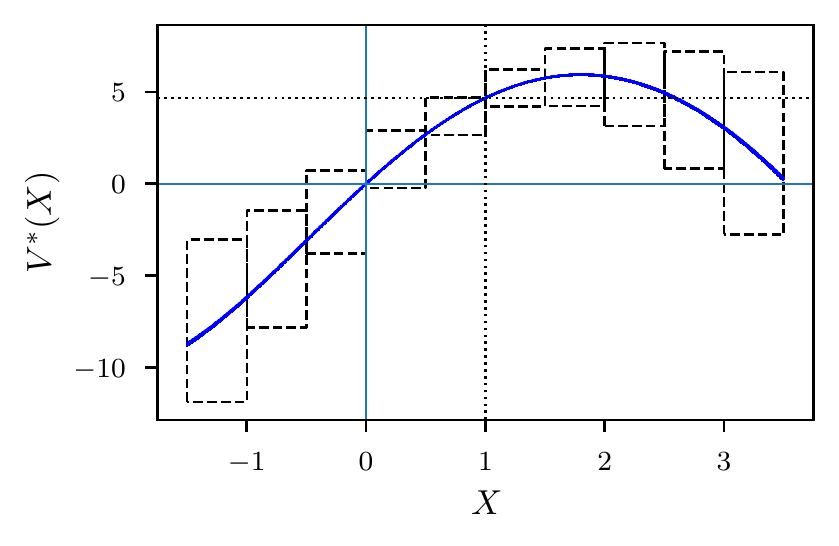}
&
\includegraphics[width=7.5cm]{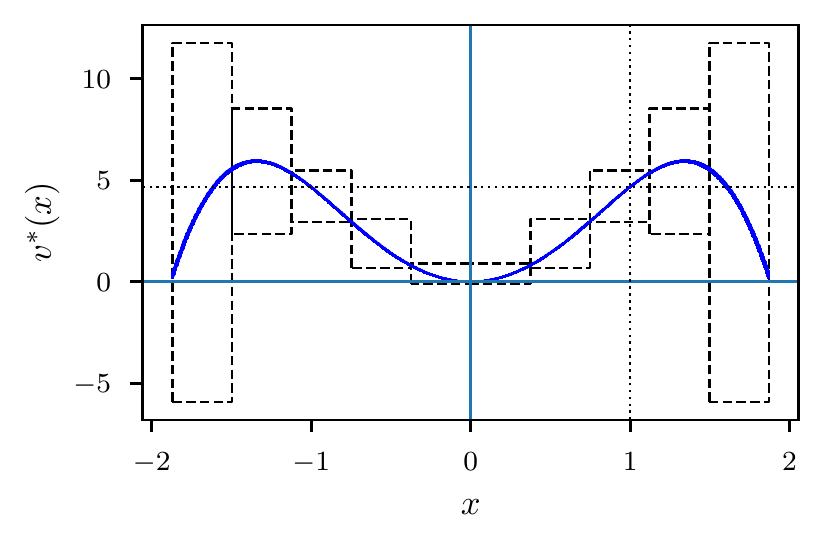}
\\
(a) & (b)\\
\includegraphics[width=7.5cm]{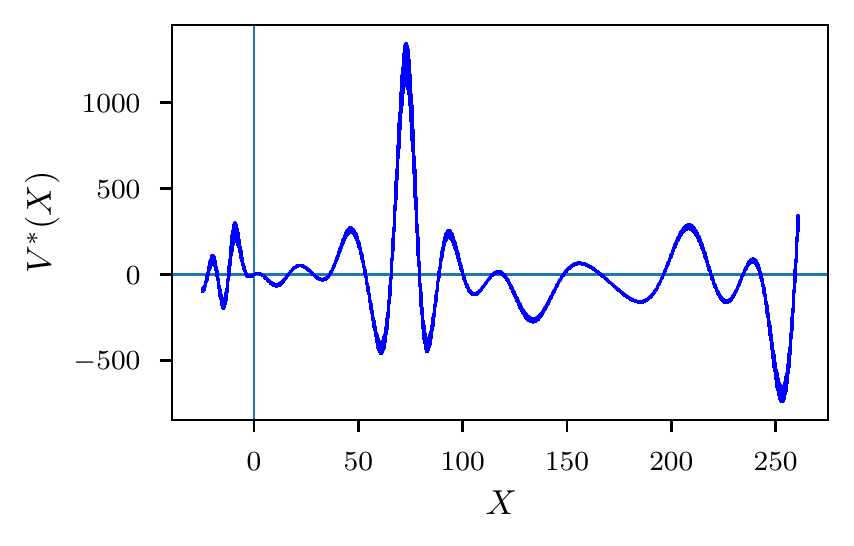}
&
\includegraphics[width=7.5cm]{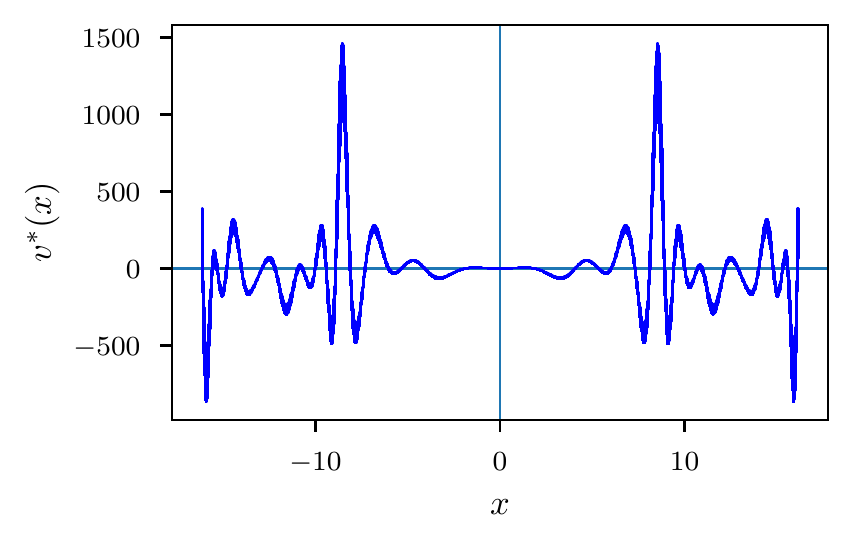}
\\
(c) & (d)
\end{tabular}
\caption{Rigorous coverings of the graphs of (a,c) the eigenfunction $V^{*}(X)$ of the linearisation $DT(G^{*})$, of $T$ at the renormalisation fixed point $G^{*}$ and (b,d) the eigenfunction $v^{*}(x)$ of $DR(g^{*})$, corresponding to eigenvalue $\delta$: (a,b) plotted for $X\in\Omega\cap\mathbb{R}$ and $x$ in the preimage (for $X\ge 0$) under $Q$, respectively, computed using a covering of the interval by $10$ (dashed lines) and $1,000$ (solid lines) subintervals, and (c,d) extended to larger domains by making use of the fixed-point equation and eigenproblem equation recursively. In (a,b) dotted lines demonstrate how the eigenvalue $\delta\simeq 4.669$ is encoded via a coordinate functional $\delta=\varphi(V^{*})$; here, $\delta=V^{*}(c)$, with $c=1$, corresponding to the constant term of the relevant power series, equation~(\ref{eqn:powerseries}), expanded with respect to the disc $\Omega=D(c,r)$; recall that $V^{*}\in\mathcal{A}(\Omega)$.\label{fig:deltaeigenfunction}}
\end{center}
\end{figure}

Compactness of $DT(G)$ implies that its spectrum  consists of $0$ together with countably-many isolated eigenvalues, each of finite multiplicity, accumulating at $0$.
We note that the spectra of $DT(G^{*})$ and $DR(g^{*})$ are related, with $\alpha^2\approx 6.264547$ and $\delta\approx 4.669201$ the only eigenvalues of $DT(G^{*})$ outside the closed unit disc, and $\alpha$ and $\delta$ the only eigenvalues of $DR(g^{*})$ outside the closed unit disc.  We recall that $\alpha$ controls universal scaling in the state variable and $\delta$ in the parameter for families of unimodal maps in the corresponding universality class.

We use a method (introduced in \cite{Bur20}), novel in the context of bounding the spectrum of derivatives of renormalisation operators, where we write an eigenvalue, $\lambda$, as a linear functional of the corresponding eigenfunction, $\lambda = \varphi(V)$, thereby expressing the eigenproblem for $(V,\lambda)$ as the following nonlinear problem for $V\in\mathcal{A}(\Omega)$:
\begin{equation}
DT(G)V-\varphi(V) V=0,\label{eqn:deltaprob}
\end{equation}
and then adapt the techniques of section~\ref{sec:renormfp} to establish that a suitably-chosen Newton-like operator $\hat{\Phi}$ for this modified problem is a contraction mapping on a ball $B(V^0,\hat{\rho})$ centred on an approximate eigenfunction $V^0$. Note that, in equation~(\ref{eqn:deltaprob}), $G$ is taken to range over the ball $B(G^0,\rho)$, proven to contain the fixed point $G^{*}$. This necessarily places a restriction on the tightness of bounds on the eigenfunction $V$ (and the corresponding eigenvalue $\lambda$) that can ultimately be achieved by this approach.
The linear functional $\varphi$ is chosen as the coordinate functional that extracts the first non-zero power series coefficient of the desired eigenfunction (when expanded with respect to $\Omega$). Our choice of norm~(\ref{eqn:powerseries},\ref{eqn:norm}) thus implies that $\hat{\rho}$ provides bounds on the eigenvalue $\lambda$ directly via~$\lambda\in[\varphi(V^0)-\hat{\rho},\ \varphi(V^0)+\hat{\rho}]$.

For the eigenfunction $V^{*}$ corresponding to the universal constant $\delta$, we were able to achieve $\hat\rho=10^{-403}$, giving $\|\hat{\Phi}(V^0)-V^0\|\le\hat\varepsilon<1.2\times 10^{-404}$, and $\|D\hat{\Phi}(V)\|\le\hat\kappa<2.8\times 10^{-100}$ for all $V\in B(V^0,\hat{\rho})$.   As in section~\ref{sec:renormfp}, care must be taken to minimise dependencies on the argument in the corresponding derivatives $D\hat{\Phi}(V)\delta V$ when bounding the action of $D\hat{\Phi}(V)$ on high-order perturbations in the computation of $\hat{\kappa}$.
The eigenfunction (with our chosen normalisation $\delta=V^{*}(1)$) is illustrated in figure~\ref{fig:deltaeigenfunction}.

As a result, we find rigorous bounds on $\delta$ that confirm $403$ significant digits
(see appendix \ref{d2summary}), with the first 20 given as:
\begin{equation}
    \delta = 4.66920\,16091\,02990\,6718\ldots
\end{equation}
We note that the same method may be employed to bound the eigenfunction corresponding to $\alpha^2$ and hence to bound $\alpha$ itself, albeit in a less direct way (and less tightly) than in section~\ref{sec:renormfp}. The results confirm previous numerical estimates.

\section{Critical scaling of added uncorrelated noise}
\label{sec:LW}

\begin{figure}[ht]
\begin{center}
\begin{tabular}{cc}
\includegraphics[width=7.5cm]{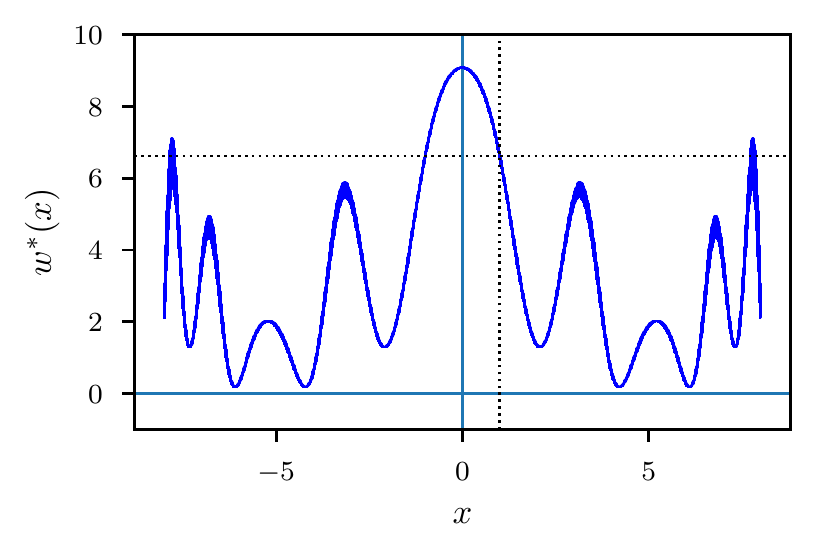}
&
\includegraphics[width=7.5cm]{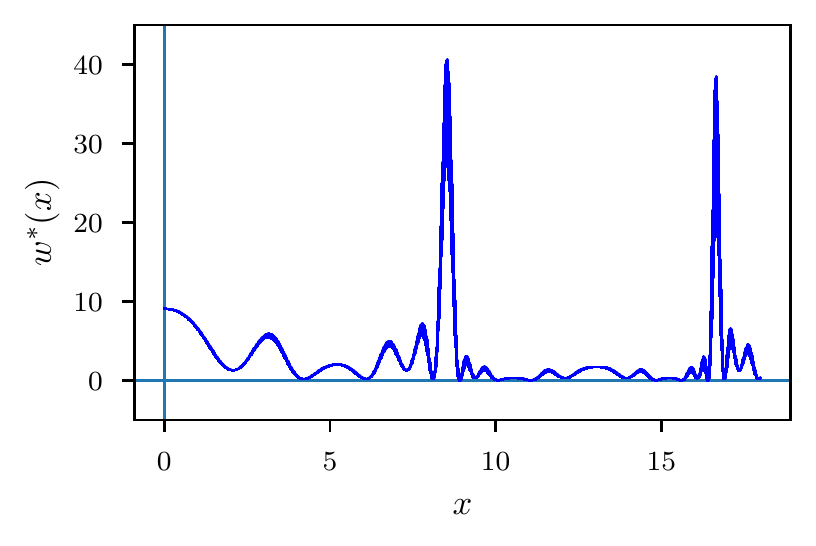}
\\
(a) & (b)
\end{tabular}
\caption{Rigorous coverings of the graph of $w^{*}(x)=W^{*}(Q(x))$ where $W^{*}$ is the eigenfunction of the operator $L$ corresponding to the eigenvalue $\gamma^2$ of largest absolute value. The function $W^{*}$ (resp. $w^{*}$) is extended to domains larger than $\Omega\cap\mathbb{R}$ (resp.its preimage for $X\ge 0$) by using the fixed-point and eigenproblem equations recursively. In (a) dotted lines demonstrate how the eigenvalue $\gamma\simeq 6.619$ is encoded via a coordinate functional $\gamma=\varphi(W^{*})$; here, $\gamma=W^{*}(c)$ with $c=1$, corresponding to the constant term of the relevant power series, equation~(\ref{eqn:powerseries}); recall that $W^{*}\in\mathcal{A}(\Omega)$.\label{fig:noiseeigenfunction}}
\end{center}
\end{figure}

We further apply the techniques of sections~\ref{sec:renormfp}--\ref{sec:DTG} to find tight rigorous bounds on the eigenfunction, $w^{*}$, and eigenvalue, $\gamma$, controlling the universal scaling of added uncorrelated noise. We must first adapt the formulation to our modified operator $T$. Following \cite{Cru81}, we consider the iteration of a prototypical one-parameter family, $f_{\mu},$ modified by the addition of independently identically distributed random variables $\xi_n$ at each iterate, to give
\begin{equation}
    x_{n+1}=F_{\mu,n}(x_n):= f_{\mu}(x_n)+\varepsilon\xi_n
\end{equation}
The relevant operator controlling universal scaling of the noise, obtained by considering the limits $\varepsilon\to 0$ and $g\to g^{*}$, is given by
\begin{equation}
    N(g)w(x)=a^{-2}\left[(g'(g(ax)))^2w(ax)+w(g(ax))\right],
\end{equation}
with eigenproblem written as
\begin{equation}
    N(g)w(x)= \lambda^2 w(x),
\end{equation}
in which the eigenvalue $\gamma^2$ of largest absolute value controls universal scaling of the overall variance of the added noise, with the corresponding eigenfunction $w(x)$ controlling the distribution in $x$.

Firstly, we note that the relevant operator, written in terms of $G$ and $X$, corresponding to our modified renormalisation operator $T$ is given by
\begin{equation}
    LW(X)=a^{-2}(G'(Q(G(a^2X)))\cdot 2G(a^2X))^2W(a^2X)+a^{-2}W(Q(G(a^2X)))
    \label{LK}
\end{equation}
where we have taken $W\circ Q=w$. The first and second terms on the right hand side of equation~(\ref{LK}) correspond to those terms in the expression for the Fr\'echet derivative $DT(G)$ that don't involve variations in $a$, namely terms~(\ref{eq:term3}) and~(\ref{eq:term2}), respectively. (We note that the terms corresponding to~(\ref{eq:term1}) and~(\ref{eq:term4}) are removed in the derivation in~\cite{Cru81} as a consequence of the independence of the $\xi_n$, and that the correlated case is dealt with in~\cite{Sch81}.)

We adapt our procedure from section~\ref{sec:DTG} for $DT(G)$, writing the eigenvalue as a linear coordinate functional of the eigenfunction to give the (nonlinear in $W\in\mathcal{A}(\Omega)$) problem
$$
LW(X)-\varphi(W)^2W=0,
$$
and use a contraction mapping argument for an appropriate Newton-like operator $\tilde{\Phi}$ on a ball $B(W^0,\tilde{\rho})$. Note, again, that in formulating the Newton-like operator corresponding to equation~(\ref{LK}) we allow $G$ to vary over the ball $B(G^0,\rho)$, proven to contain the renormalisation fixed point $G^{*}$. We thereby establish rigorous bounds on both the eigenfunction $W^{*}$ and eigenvalue $\gamma=\varphi(W^{*})$.

Using a ball of radius $\tilde\rho=10^{-401}$, we gain bounds $\|\tilde{\Phi}(W^0)-W^0\|\le\tilde\varepsilon<8.9\times 10^{-402}$, and $\|D\tilde{\Phi}(W)\|\le\tilde\kappa<7.4\times 10^{-101}$ for all $W\in B(W^0,\tilde{\rho})$. The eigenfunction (with our chosen normalisation) is illustrated in figure~\ref{fig:noiseeigenfunction}.
This yields $401$ correct significant digits for $\gamma$ 
(see appendix \ref{d2summary}), with the first 20 given as
\begin{equation}
    \gamma = 6.61903\,65108\,17928\,0453...
\end{equation}

\clearpage
\appendix
 \section{Appendix}
 \subsection{Tight rigorous bounds on universal constants}
 \label{d2summary}
 Below we state the digits proven correct of $a=\alpha^{-1}=G^{*}(1)$ ($409$ digits), $\alpha$ ($408$ digits), $\delta=\varphi(V^{*})=V^{*}(1)$ ($403$ digits), and $\gamma=\varphi(W^{*})=W^{*}(1)$ ($401$ digits) obtained from a proof with truncation degree $N=640$ for the polynomial part of the space for $G^*$, $V^*$, and $W^*$ (corresponding to degree $2N = 1280$ for $g^*$, $v^*$, and $w^*$, respectively). The  numbers listed here are found by computing safely rounded decimal interval bounds, containing the corresponding dyadic rational computer-representable bounds, on the relevant constants, and quoting only those digits guaranteed accurate as a result.

\begin{center}
\begin{tabular}{rlllll}
$a={}$
-0.
&3995352805 &2313448985 &7580468633 &6937194335 &4428046695\\ &2727517073 &0449124380 &1660883804 &2981844594 &8741812667\\
&6179406484 &6838366714 &0945404846 &1643643736 &0947557018\\
&4545976789 &4023268702 &2548579773 &5028209746 &4775103925\\
&5797877507 &3697474932 &3269755137 &3492308212 &2088541722\\
&2413083309 &4802739189 &0574703944 &6460416066 &9938415778\\
&2298900077 &7299013544 &2121397192 &4552385259 &4449033723\\
&7697553775 &0905488329 &7544336726 &9368114050 &5788840461\\
&79344018\ldots\\

$\alpha={}$
-2.
&5029078750 &9589282228 &3902873218 &2157863812 &7137672714\\ &9977336192 &0567792354 &6317959020 &6703299649 &7464338341\\ &2959523186 &9995854723 &9421823777 &8544517927 &2863314993\\ &3725781121 &6359487950 &3744781260 &9973805986 &7123971173\\ &7328927665 &4044010306 &6983138346 &0009413932 &2364490657\\ &8899512205 &8431725078 &7337746308 &7853424285 &3519885875\\ &0004235824 &6918740820 &4281700901 &7148230518 &2162161941\\ &3199856066 &1293827426 &4970984408 &4470100805 &4549677936\\ &7608881\ldots\\

$\delta={}$
+4.
&6692016091 &0299067185 &3203820466 &2016172581 &8557747576\\ &8632745651 &3430041343 &3021131473 &7138689744 &0239480138\\ &1716598485 &5189815134 &4086271420 &2793252231 &2442988890\\ &8908599449 &3546323671 &3411532481 &7142199474 &5564436582\\ &3793202009 &5610583305 &7545861765 &2222070385 &4106467494\\ &9428498145 &3391726200 &5687556659 &5233987560 &3825637225\\ &6480040951 &0712838906 &1184470277 &5854285419 &8011134401\\ &7500242858 &5382498335 &7155220522 &3608725029 &1678860362\\ &67\ldots\\

$\gamma={}$
+6.
&6190365108 &1792804532 &3808905147 &4666014364 &4298809101\\ &1980889058 &1539120755 &2294388390 &1250134543 &0103013791\\ &0116621507 &6680991461 &7111062123 &4676765967 &2263346641\\ &5349015651 &5469980646 &2621251411 &3242709973 &9377082075\\ &2957874751 &6962711711 &6928533607 &9067798211 &8951469414\\ &0224500385 &6708624240 &5473933494 &6093414214 &2285269246\\ &7145643730 &8826640353 &2825154865 &5386124267 &3930589439\\ &2213420488 &3953151516 &3766198410 &1165280871 &0270346725\\ &\ldots
\end{tabular}
\end{center}
\subsection{Computational Considerations}
The computations were performed via two independent implementations of the rigorous framework: (1) written in the high-performance language Julia, making use of IEEE754-2008 compliant binary multi-precision arithmetic with rigorous directed rounding, and (2) written in the language Python, using IEEE754-2008 compliant decimal multi-precision arithmetic with rigorous directed rounding. In both cases, directed rounding modes are respected at the process level, and multiprocessing rather than threads was therefore used to ensure safety during parallel computations.

\end{document}